\documentclass{scrartcl}

\KOMAoptions{paper=a4}
\KOMAoption{fontsize}{10pt}

\usepackage[utf8]{inputenc} 
\usepackage[T1]{fontenc} 
\usepackage[ngerman,english]{babel} 
\usepackage{hyphenat}

\usepackage{lmodern}
\usepackage[scaled=0.9]{helvet}
\usepackage{courier}

\usepackage{todonotes} 

\usepackage{enumitem}
\setlist[enumerate]{label*=(\alph*),ref=(\alph*)}

\usepackage{graphicx} 
\usepackage{float} 
\graphicspath{{pic/}} 

\usepackage{amsmath} 
\usepackage{amssymb} 
\usepackage{amscd} 
\usepackage{amsthm}

\usepackage{csquotes} 
\usepackage[citestyle=alphabetic,
            bibstyle=alphabetic,
						sorting=nyt,
						maxnames=4,
						backend=biber,
						isbn=false,
						doi=false,
						clearlang=true]
						{biblatex}
\bibliography{biblio}
\DeclareFieldFormat*{title}{\textit{#1}}
\DeclareFieldFormat{journaltitle}{#1}
\DeclareFieldFormat{booktitle}{#1}
\DeclareRedundantLanguages{english,English}{english,german,ngerman,french}

\usepackage{hyperref} 
\usepackage{color} 
\hypersetup{citecolor=blue, urlcolor=blue, linkcolor=red, colorlinks=true, breaklinks=true} 

\DeclareMathAlphabet{\mathds}{U}{dsrom}{m}{n}
\newcommand{\CC}{{\mathds C}}

\newcommand{\PP}{{\mathds P}}

\newcommand{\RR}{{\mathds R}}

\newcommand{\ZZ}{{\mathds Z}}

\DeclareMathOperator {\Star}{Star}

\DeclareMathOperator {\divisor}{div}

\DeclareMathOperator {\id}{id}

\newcommand{\Fans}{\text{Fans}}
\DeclareMathOperator{\lindim}{lindim}
\DeclareMathOperator{\spldim}{spldim}
\DeclareMathOperator{\LinSp}{LinSp}
\DeclareMathOperator{\Rec}{Rec}
\DeclareMathOperator{\RelInt}{RelInt}
\DeclareMathOperator{\bsim}{{\stackrel{b}{\sim}}}
\DeclareMathOperator{\Rsim}{{\stackrel{\RR}{\sim}}}
\newcommand{\dfn}[1]{\emph{#1}}
\newcommand{\sss}{{\mathfrak s}}

\newcommand{\TT}{{\mathds T}}
\newcommand{\XXX}{{\mathcal X}}

\newcommand{\FFF}{{\mathcal F}}

\newtheorem {theorem}{Theorem}[section]
\newtheorem {proposition}[theorem]{Proposition}
\newtheorem {lemma}[theorem]{Lemma}

\newtheorem {corollary}[theorem]{Corollary}

\theoremstyle {definition}
\newtheorem {definition}[theorem]{Definition}
\newtheorem {example}[theorem]{Example}

\theoremstyle {remark}
\newtheorem {remark}[theorem]{Remark}

\begin {document}

\title {On rational equivalence in tropical geometry}
\author {Lars Allermann, Simon Hampe, Johannes Rau}
\date{} 

\maketitle

\begin{abstract}
  This article discusses the concept of rational equivalence in tropical
	geometry (and replaces the older and imperfect version \cite{AR2}). 
	We give the basic definitions in the context of tropical varieties
	without boundary points and prove some basic properties. 
	We then compute the ``bounded'' Chow groups of $\RR^n$ by showing that they are isomorphic
	to the group of fan cycles. The main step in the proof is of independent interest:
	We show that every tropical cycle in $\RR^n$ is a sum of (translated) fan cycles. This also 
	proves that the intersection ring of tropical cycles is generated in codimension 1 (by hypersurfaces).
\end{abstract}


\section{Introduction}

The concept of rational equivalence plays a fundamental role in algebraic geometry
and therefore it is natural to study analogue notions in tropical geometry.
This has been done quite extensively in the case of divisors on a curve
(e.g.\ \cite{BN, GK, MZ2, CDPR, HMY}), whereas in higher dimensions
there are relatively few instances where rational equivalence is mentioned explicitly
(cf.\ \cite{Mi3, AR, MR2}, for example). 

This paper is devoted to the basic definitions and properties of rational equivalence for 
tropical varieties. We stick to non-compact tropical varieties without ``boundary'' points here
and study usual as well as ``bounded'' rational equivalence 
(the latter using bounded rational functions). 
We prove some basic properties (in particular the compatibility with the constructions from \cite{AR})
and show that bounded rational equivalence can also be expressed in terms
of families of cycles over $\RR$. 

We then turn to the case of cycles in $\RR^n$ and show that two cycles are bounded rationally equivalent
if and only if they are numerically equivalent if and only if they have the same recession fan cycle.
It follows that the bounded Chow group of $\RR^n$ is isomorphic to the group of fan cycles in $\RR^n$.
The main step is to prove that a tropical cycle is rationally equivalent to its recession fan cycle.
We deduce this by proving another statement of independent interest: We show that every tropical cycle
in $\RR^n$ can be decomposed into a sum of (translated) fan cycles. This also proves the fact
that every such tropical cycle can be written as a sum of intersection products of hypersurfaces.
In other words, hypersurfaces $V(f)$ with $f$ a tropical polynomial generate the ring of tropical
cycles $Z_*(\RR^n)$.

There exists an older and imperfect version of this paper by the first and third author on arXiv
(cf.\ \cite{AR2}). 
Our main motivation for this new version was to replace the proof of the ``main step'' mentioned above 
(i.e.\ rational equivalence of a cycle and its recession fan cycle) by a simpler and more transparent argument.
To us, the cleanest way in order to update the old paper seemed to be to replace it completely
and therefore to include the old material in this new version. 
In doing so, we also updated the terminology slightly. 
Rational equivalence in the old paper is now called \emph{bounded rational equivalence} (as it is 
generated by bounded rational functions). We added the concept of usual \emph{rational equivalence}
(generated by (arbitrary) rational functions) and \emph{rational equivalence over $\RR$} (generated
by families over $\RR$). 

The authors would like to thank
Andreas Gathmann and Hannah Markwig 
for useful comments and discussions.

\section{Preliminaries}

This article is, to some extent, a continuation of \cite{AR} and
we mostly stick to the definitions and notations introduced there. 
However, for the reader's convenience we start by recalling the most important 
terminology. For more details, we kindly refer the reader to the cited work.

\subsection {Cycles} \label{cycles}


A \emph{tropical polyhedral complex $\XXX$} is a balanced (weighted, pure-dimensional, rational, finite) polyhedral complex 
in $\RR^n$ (with underlying lattice $\ZZ^n$). 
The top-dimensional polyhedra
in $\XXX$ are called \emph{facets}, the codimension one polyhedra are
called \emph{ridges}. \emph{Balanced} means that for each ridge $\tau \in \XXX$ the following 
\emph{balancing condition at $\tau$} is satisfied: The weighted
sum of the primitive vectors of the 
facets $\sigma$ around $\tau$ is zero, i.e.\
\begin{equation*} 
  \sum_{\substack{\sigma \in \XXX^{(\dim(\XXX))} \\ \tau < \sigma}} 
    \omega(\sigma) v_{\sigma / \tau} = 0.  
\end{equation*}  
Here, $\omega(\sigma)$ denotes the weight of the facet, and
$v_{\sigma / \tau}$ is the primitive integer generator of the 
ray obtained from projecting $\sigma$ to $\RR^n / V_\tau$, where 
$V_\tau$ denotes the linear vector space spanned by $\tau$.
The \emph{support of $\XXX$}, denoted by $|\XXX|$, is the union of all facets
in $\XXX$ with non-zero weight. 


Two tropical polyhedral complexes are called equivalent if the they admit a 
common refinement and if the induced weights are the same. 
A \emph{tropical cycle $X$} is an equivalence class of tropical polyhedral complexes. 
A representative $\XXX$ of $X$ is called a polyhedral structure for $X$. 
Obviously, the support of $X$ is well-defined and we
often denote it by the same letter $X$.
Consistently with this abuse of notation, we may think of a tropical cycle $X$ as a polyhedral set 
with weights $\omega_X(p)$ for generic points $p \in X^{\text{gen}}$ 
such that (after choosing a polyhedral structure) the balancing condition is satisfied.  
A tropical cycle $F$ supported on a fan (i.e.\ a union of cones with vertex at $0$)
is called a \emph{fan cycle}.




%
%

\subsection{The divisor of a rational function} 

\label{principaldivisor}


A map $\varphi : \RR^n \supseteq S \to \RR^m$ is called \dfn{integer affine}
if there exist $A \in \text{Mat}(m \times n,\ZZ)$ and $a \in \RR^m$ such that 
for all $p \in S$
\begin{equation*} 
  \varphi(p) = A p + a.
\end{equation*}
A \emph{(non-zero) rational function} on a tropical cycle 
$X$ is a continuous function $\varphi : X \rightarrow \RR$ 
that is integer affine on each cell of a suitable 
polyhedral structure $\XXX$ of $X$.
The \emph{divisor of $\varphi$}, denoted by
$\divisor(\varphi) = \varphi \cdot X$, is given by the
weighted subcomplex $\varphi \cdot \XXX$ of $\XXX$ constructed in \cite[\nopp 3.3]{AR}.
It is supported on the codimension one skeleton of $\XXX$ and 
contains each ridge $\tau \in \XXX$ (now a facet of $\varphi \cdot \XXX$)
with weight
\begin{equation} \label{weightdivisor}
  \omega_{\varphi \cdot X}(\tau) =
    \sum_{\substack{\sigma \in X^{(\dim X)} \\ \tau < \sigma}} 
    \omega(\sigma) \varphi_\sigma(\tilde{v}_{\sigma / \tau}) 
    - \varphi_\tau\Big( \sum_{\substack{\sigma \in X^{(\dim X)} \\ \tau < \sigma}} 
    \omega(\sigma) \tilde{v}_{\sigma / \tau}\Big).
\end{equation}
Here $\varphi_\sigma : V_\sigma \rightarrow \RR$ denotes the 
linear part of the affine function
$\varphi|_\sigma$, and $\tilde{v}_{\sigma / \tau}$ is an arbitrary representative  in $\RR^n$ of
$v_{\sigma / \tau} \in \RR^n/V_\tau$.
It was shown in \cite[\nopp 3.7]{AR} that these weights satisfy the 
balancing condition, hence $\divisor(\varphi)$ is well-defined 
tropical tropical subcycle of $X$ of codimension one. 
Note also that $\divisor(\varphi)$ agrees with the intersection
of the \emph{balanced graph of $\varphi$}
with $X \times \{-\infty\}$.
The balanced graph of $\varphi$ is obtained from the usual graph of $\varphi$
(not balanced, in general) by adding cells in the $(0,\ldots, 0,-1)$-direction in order
to make it balanced. 
In this sense, $\divisor(\varphi)$ can be regarded
as the divisor of zeros and poles (if negative weights show up) of
$\varphi$.

\subsection{Morphisms and projection formula} \label{projection}


Given two cycles $X \subseteq \RR^n$ and $Y \subseteq \RR^m$, a integer
affine map $f : X \rightarrow Y$ is called a \dfn{morphism of cycles}.
Given such a morphism, we can \dfn{pull back} a rational function $\varphi$ on $Y$ 
to a rational function $f^*(\varphi) = \varphi \circ f$ on $X$. 
Furthermore, we can \dfn{push forward} a subcycle 
$Z$ of $X$ to a subcycle $f_*(Z)$ of $Y$. This is due to
\cite[2.24 and 2.25]{GKM} in the case of fans and can be generalized
to complexes (see \cite[\nopp 7.3]{AR}). 
The push forward 
$f_*(Z)$ is supported on the image $f(|Z|)$ and
(for sufficiently fine polyhedral structures) the weights of $f_*(Z)$ 
are given by 
\begin{equation*} 
  \omega_{f_*(Z)}(\sigma') =
    \sum_{\substack{\sigma \text{ facet of } Z \\ f(\sigma) = \sigma'}}
    |\Lambda_{\sigma'} / f(\Lambda_\sigma)| \cdot \omega_Z(\sigma), 
\end{equation*}
where $\sigma'$ is a cell of $Y$ of dimension $dim(Z)$.
Here, $\Lambda_\sigma := V_\sigma \cap \ZZ^n$ denotes the 
sublattice of $\ZZ^n$ spanned by $\sigma$ (analogously for $\Lambda_{\sigma'}$). 
It follows that $\dim(f_*(Z)) = \dim(Z)$ if $f_*(Z) \neq 0$.

The \emph{projection
formula} (see \cite[\nopp 4.8]{AR}) connects all the above constructions via
\begin{equation} \label{eq:projectionformula}
  f_*(f^*(\varphi) \cdot Z) = \varphi \cdot f_*(Z).
\end{equation}

\subsection{Intersection product of two cycles}
\label{productofcycles}

Another feature of tropical intersection theory is
that for any two cycles $X, Y$ in $\RR^n$ we can perform
the ``stable intersection'' $X \cdot Y$ which is again a
well-defined cycle in $\RR^n$ (not just a cycle class modulo rational equivalence).
The codimension of $X \cdot Y$ is always equal to the sum of codimensions of $X$ and $Y$, 
regardless of the dimension of the set-theoretic intersection $X \cap Y$.
The definition given in \cite[\nopp 9.3]{AR} 
is based on intersecting the cartesian product $X \times Y \subseteq \RR^n \times \RR^n$ with
the diagonal described by the rational functions 
$\max\{x_1,y_1\}, \ldots, \max\{x_n,y_n\}$, i.e.\
\begin{equation*} 
  X \cdot Y := \pi_*\big(\max\{x_1,y_1\} \cdots \max\{x_n,y_n\} 
    \cdot (X \times Y)\big).
\end{equation*}
Here, $x_i, y_i$ are the coordinates of the first resp.\ second factor of $\RR^n$
and $\pi$ is any of the two projections.
This intersection product turns $Z_*(\RR^n)$ into a graded commutative ring, and satisfies 
$(\varphi \cdot X) \cdot Y = \varphi \cdot (X \cdot Y)$, where $\varphi$ is a rational function
on $X$. Moreover, $\RR^n$ (considered as a cycle) is the identity element. 

Note that our definition is a way of formalizing the concept stable intersection based on moving the cycles
slightly as proposed in \cites{RST,Mi3} (details can be found e.g.\ in \cite{MR};
equivalence was proven in \cites{Ka,Ra}).

\section{Rational equivalence}

As discussed in \cite[\nopp 8.6]{AR}, the definition of rational equivalence 
given there is not compatible with push-forwards of cycles. 
The following definition is more flexible and resolves this problem.
Moreover, we show in proposition \ref{pro:seconddefinition} that
this definition is consistent with the approach of using families 
over $\RR$.


\begin{definition} \label{ratequ1}
  Let $X$ be a cycle and let $Z$ be a subcycle. We call
  $Z$ \emph{bounded rationally equivalent to zero on $X$}
  if there exists a morphism $f : Y \rightarrow X$ and a
  \emph{bounded} rational function $\phi$ on $Y$ such that
	\begin{equation*} 
    f_*(\phi \cdot Y) = Z.
	\end{equation*}
	Note that in this case $\dim(Y) = \dim(Z) + 1$. 
  Two subcycles $Z, Z'$ of $C$ are called 
	\dfn{rationally equivalent}, denoted by $Z \sim Z'$, if $Z-Z'$ is rationally
  equivalent to zero.
	
	Furthermore, if the function $\phi$ from above can chosen to be bounded,
	we call $Z$ \dfn{bounded rationally equivalent to zero}. 
	The corresponding equivalence relation is called \dfn{bounded rational equivalence}
	and is denoted by $\bsim$.
\end{definition}

Obviously, both $\sim$ and $\bsim$ are additive equivalence relations.
It also clear that $Z \bsim Z'$ implies $Z \sim Z'$, hence
bounded rational equivalence is the stronger relation.
Even though $\sim$ is more natural from the classical point of view,
the main interest in this paper will be on bounded rational equivalence.
Some explanations regarding this are collected in remark \ref{rem:noboundary}.
(Bounded) rational equivalence satisfies the following properties.

\begin{proposition} \label{RationalEquivalenceProperties}
  Let $Z$ be a cycle in $X$ (bounded) rationally equivalent to zero.
  Then the following holds:
  \begin{enumerate}
    \item
      Let $X'$ be another cycle. Then
      $Z \times X' \subseteq X \times X'$ is also (bounded) rationally equivalent to zero.
    \item
      Let $\varphi$ be a rational function on $X$.
      Then $\varphi \cdot Z$ is also (bounded) rationally equivalent
      to zero.
    \item
      Let $g : X \rightarrow \widetilde{X}$ be a morphism. Then
      $g_*(Z)$ is also (bounded) rationally equivalent to zero.
    \item
      Assume $X = \RR^n$ and let $Z'$ be another
      cycle in $\RR^n$. Then $Z \cdot Z'$
      is also (bounded) rationally equivalent to zero. 
		\item
      If $Z \bsim 0$ and $\dim(Z)=0$, then $\deg(Z) = 0$. 
			Here, as usual, $\deg(\sum m_i P_i) := \sum m_i$ denotes the sum 
			of coefficients.
  \end{enumerate}
\end{proposition}

\begin{proof}
  Let $f : Y \rightarrow X$ be a morphism and $\phi$
  a (bounded) function on $Y$ such that $f_*(\phi \cdot Y) = Z$.
  Then $f \times \id : Y \times X' \rightarrow X \times X'$ shows (a).
	Restricting $f$ to $f^*(\varphi) \cdot Y$ and using 
	the projection formula \eqref{eq:projectionformula} shows (b).
	Composing $f$ with $g$ shows (c). 
  For (d) we just have to recall that $Z \cdot Z'$ is computed 
  by
	\begin{equation*} 
    \pi_*(\max\{x_1,y_1\} \cdots \max\{x_n,y_n\} \cdot (Z \times Z')).
	\end{equation*}
	Thus (d) follows from (a) -- (c). 
	We are left with (e), the only case where the stronger concept
	of bounded rational equivalence is needed. 
	In this case, $Y$ must be one-dimensional and
  we can apply \cite[\nopp 8.3]{AR}, which shows that
  the degree of $\phi \cdot Y$, when $\phi$ is bounded, is zero. Pushing forward 
	preserves degree, and hence the statement follows.
\end{proof}

An easy example of bounded rationally equivalent cycles are translations.
Given a vector $\vec{v} \in \RR^n$, in the following
$X + \vec{v}$ will always denote the translation of a cycle
$X$ by the vector $\vec{v}$. This should be distinguished
from the usual sum of cycles $X + Y$ given taking unions and adding weights.

\begin{proposition} \label{TranslationEquivalent}
  Let $X$ be a cycle in $\RR^n$ and let $X + \vec{v}$ denote the
  translation of $X$ by an arbitrary
  vector $\vec{v} \in \RR^n$. Then 
	\begin{equation*} 
    X \bsim X + \vec{v}.
	\end{equation*}
\end{proposition}

\begin{proof}
  Consider the
  cycle $X \times \RR$ in $\RR^n \times \RR$ and the 
  morphism
  \begin{eqnarray*}
    f : \RR^n \times \RR & \rightarrow & \RR^n, \\
        (x, t)       & \mapsto     & x + t  \vec{e}_i,
  \end{eqnarray*}
  where $\vec{e}_i$ is the $i$-th unit vector in $\RR^n$.
  For $\mu \in \RR_\geq$ let $\phi_\mu$ be the bounded
  function
	\begin{equation*} 
    \phi_\mu(x,t) = \begin{cases}
                    0 & t \leq 0 \\
                    t & 0 \leq t \leq \mu \\
                    \mu & t \geq \mu.
                  \end{cases}
	\end{equation*}
  Then we can compute 
	\begin{equation*} 
		\begin{split}
  		f_*(\phi_\mu \cdot (X \times \RR)) &= f_*(C \times \{0\} - C \times \{\mu\}) \\
	  	                                   &= C - (C + \mu \vec{e}_i).
		\end{split}
	\end{equation*}
	Applying this to each coordinate step by step, we obtain $X \bsim X + \vec{v}$.
\end{proof}

In algebraic geometry, instead of using the 
divisors of zeros and poles of rational functions, 
one may define rational equivalence 
by considering (flat) families
of cycles over $\PP^1$. Indeed, two cycles are rationally equivalent if
they both appear as fibers of some family $F$. 
Tropically, we can do the same (cf.\ \cite[\nopp 4.6]{Mi3}).


\begin{definition} \label{ratequ2}   
  Let $X$ be a cycle and consider a subcycle $F \subseteq X \times \RR$.
	For each point $p \in \RR$, we define the \dfn{fiber} of $F$ at $p$ by
	\begin{equation*} 
	  F_p := \varphi_p \cdot F \subseteq X \times \{p\} \cong X,
	\end{equation*}
	where $\varphi_p$ is the pull-back of $\max\{x,p\}$ along $X \times \RR \to \RR$.
	We think of $F_p$ as a subcycle of $X$.
	The equivalence relation generated by setting 
	$F_p \Rsim F_q$ for two fibers $F_p, F_q$ of the same $F$
	is called \dfn{rational eqivalencce over $\RR$}, denoted by $\Rsim$.
\end{definition}

\begin{proposition} \label{pro:seconddefinition}
  Rational equivalence over $\RR$ (as defined in \ref{ratequ2}) agrees with
	bounded rational equivalence (as defined in \ref{ratequ1}).
\end{proposition}

\begin{proof}
  We first show that, given $F \subseteq X \times \RR$ and $p,q \in \RR$,
	any two fibers $F_p$ and $F_q$ satisfy $F_p \bsim F_q$. 
	To see this, let $\varphi$ be a rational function on $\RR$
	with divisor $\varphi \cdot \RR = \sum m_i p_i$. 
	Let furthermore denote $\pi: F \to X \times \RR \to \RR$ the second projection.
	Pulling back $\varphi$ to $F$, we obtain
	\begin{equation} \label{eq:sumoffibers}
	  \pi^* \varphi \cdot F = \sum m_i F_{p_i}.
	\end{equation}
	This follows from the fact that the divisor construction is local (cf.\
	\cite[\nopp 1.1]{Ra}), linear and invariant under change 
	by an affine function (cf.\ \cite[\nopp 3.6]{AR}). 
	We apply equation \eqref{eq:sumoffibers} to the function
	\begin{equation*} 
	  \varphi = \max\{x,p\} - \max\{x,q\},
	\end{equation*}
	which is obviously bounded.
	We obtain
	\begin{equation*} 
	  \pi_*(\varphi \cdot F) = F_p - F_q \in Z_*(X),
	\end{equation*}
	which proves $F_p \bsim F_q$. 
	
	Let now $f : Y \to X$ be some morphism of cycles and $\varphi$ a bounded rational
	function on $Y$. We have to show $f_*(\varphi \cdot Y) \Rsim 0$. 
	In order to construct a suitable $F$, let us first consider 
	the balanced graph of $\varphi$ mentioned in \ref{principaldivisor}. 
	It is obtained from the usual graph of $\varphi$ in $Y \times \RR$
	by adding facets directed downwards in such a way that the constructed polyhedral complex satisfies 
	the balancing condition (cf.\ \cite[\nopp 3.3]{AR}). 
	Let us denote this subcycle of $Y \times \RR$ by $\Gamma$. 
	As $\varphi$ is bounded from above, we may choose $p \in \RR$ close to $+\infty$ such that
	$\Gamma_p = [\emptyset] = 0$. But $\varphi$ is also bounded from below. Hence, choosing
	$q \in \RR$ close to $-\infty$ we will only intersect the ``additional'' facets of $\Gamma$
	and therefore $\Gamma_q = \varphi \cdot Y$. 
	Let us now consider the map $f \times \id : Y \times \RR \to X \times \RR$ and set
	\begin{equation*} 
	  F := (f \times \id)_*(\Gamma) \subseteq C \times \RR. 
	\end{equation*}
	Using the projection formula, we easily see that $F_x = f_*(\Gamma_x)$ for all $x \in \RR$.
	In particular, $F_p = 0$ and $F_q = f_*(\varphi \cdot Y)$. Thus $f_*(\varphi \cdot Y) \Rsim 0$
	and we are done.
\end{proof}

In the following we will abandon the terminology ``over $\RR$'' and notation $\Rsim$
in favor of ``bounded'' and $\bsim$.

\begin{definition}
  \label{def:chowgroup}
  The \dfn{(bounded) Chow group} of $X$ is defined to be the group of tropical
	subcycles of $X$ modulo (bounded) rational equivalence, denoted by
	\begin{equation*} 
		A_*(X) := Z_*(X) / \sim \text{ and } A^\text{b}_*(X) := Z_*(X) / \bsim.
	\end{equation*}
\end{definition}

\begin{remark} \label{rem:noboundary}
  A few remarks regarding our definitions of rational equivalence might be helpful at this point.
	Note that in this paper we only work with spaces
	which do not contain ``boundary points''
	(e.g.\ the points at infinity in $\TT\PP^1 = \RR \cup \{\pm \infty\}$)
	as introduced for example in \cite{Mi3}.
	We refer to the book in progress \cite{MR} for the corresponding theory
	in this more general setting. In particular, definition 
	\ref{ratequ2} can be changed to allow families over	$\TT\PP^1$ and not just $\RR$, 
	in which case we recover rational equivalence $\sim$ (with unbounded functions). 
	Hence this equivalence relation is the canonical choice from the classical point of view. 
	However, when working with non-compact spaces, cycles can often be moved off ``to infinity''
	and hence the corresponding Chow groups contain relatively little information.
	For example, we will show $A_*(\RR^n) = 0$ (cf.\ \ref{generatedbyhypersurfaces}), 
	in analogy with the classical statement
	$A_*((\CC^*)^n) = 0$.
	
	In contrast, bounded rational equivalence in essence prohibits to move cycles to infinity
	and therefore provides richer Chow groups also in the non-compact case. 
	The main idea is that two cycles are bounded rationally equivalent in $X$ if
	and only if they are rationally equivalent in any (toric) compactification 
	$\overline{X}$ of $X$.
	For example, in the case of $X = \RR^n$ our main result \ref{ComparisonOfRationalNumericalDegree}
	together with \cite[\nopp 4.2]{FS} shows that $A^\text{b}_*(\RR^n)$ can be described as the direct limit
	of the Chow groups of all toric varieties compactifying $(\CC^*)^n$. 
	%
  %
  %
\end{remark}

\section{Numerical Equivalence}

Let us now compare bounded rational equivalence to numerical equivalence.

\begin{definition}
  Let $X$ be a cycle in $\RR^n$ of codimension $k$. Then
  we define $d_X$ to be the map
	\begin{equation*} 
		\begin{split}
    d_X: Z_k(\RR^n)      & \rightarrow  \ZZ, \\
    Z           & \mapsto      \deg(X \cdot Z).
		\end{split}
	\end{equation*}
	We call two cycles $C, D$ \emph{numerically equivalent} if the
	two functions $d_C$ and $d_D$ coincide.
\end{definition}

Note that Lemma \ref{RationalEquivalenceProperties} implies that 
bounded rationally equivalent cycles are also numerically equivalent. 
In Theorem \ref{ComparisonOfRationalNumericalDegree} we will also
prove the converse. In this section, our goal is to show
that two bounded rationally (resp.\ numerically) equivalent \emph{fan} cycles
have to be equal. 

\begin{proposition} \label{EquivalentCyclesAreEqual}
  Let $F_1$ and $F_2$ be fan cycles in $\RR^n$.
	If $F_1 \bsim
  F_2$ or $d_{F_1} = d_{F_2}$, then $F_1$ and $F_2$ are equal.
\end{proposition}

We need the following technical result.

\begin{lemma} \label{SimplicialCompletion}
  Let $F$ be a $d$-dimensional fan cycle in $\RR^n$.
	Then there exists a complete simplicial rational fan $\Theta$ in $\RR^n$
	such that $F$ can be represented by a tropical fan $\FFF$ which is a subfan of $\Theta$
	 (i.e.\ each cone of $\FFF$ is a cone of $\Theta$).
\end{lemma}

\begin{proof}
  We start with some fan $\FFF_0 = \{\sigma_1,\ldots,\sigma_N\}$ 
	representing $F$.
	Each cone $\sigma_i$ is described by certain integer linear inequalities, say 
	\begin{equation*} 
	  \sigma_i = \{x \in \RR^n : \langle f_1^i,x \rangle \geq 0,\ldots, \langle f_{k_i}^i,x \rangle \geq 0 \},
	\end{equation*}
	with $f^i_j \in \ZZ^n$. Let $H_{f^i_j}$ be the fan consisting of the two halfspaces and
	the hyperplane defined by $f^i_j$, i.e.\
	\begin{equation*} 
	  H_{f^i_j} := \left\{ \{x:\langle f_j^i,x \rangle\geq 0\}, \{x:\langle f_j^i,x \rangle=0\}, \{x:\langle f_j^i,x \rangle\leq 0\}\right\}.
	\end{equation*}
	Consider the “intersection” of all these fans, 
	$$
	  \Theta' := \bigcap_{i=1}^N \bigcap_{j=1}^{k_i} H_{f^i_j}
	$$
  as described in \cite[2.5(e)]{GKM}. In other words,
	$\Theta'$ is the complete fan in $\RR^n$ containing any cone
	which can be described by some collection of inequalities of the
	form $\pm f^i_j(x) \geq 0$. 
	By construction, $F$ can be represented by a subfan $\FFF'$ of $\Theta'$. 
	By subdividing $\Theta'$ further, we can construct 
	a simplicial fan $\Theta$ (cf.\ \cite[48]{Fu2}).
	As $\Theta$ is a refinement of $\Theta'$, $F$ can still be represented
	by a subfan of $\Theta$ (namely by $\FFF := \FFF' \cap \Theta$) and we are done.
\end{proof}

\begin{proof}[Proof of \ref{EquivalentCyclesAreEqual}]
  As mentioned before, note that $F_1 \bsim F_2$ implies $d_{F_1}=d_{F_2}$ by 
	proposition \ref{RationalEquivalenceProperties} (d) and (e).
  Hence it suffices to show the following: If $F$ is a tropical cycle with $d_F = 0$, then $F = 0$. 
  We prove this by induction on $d := \dim(F)$. For $d=0$ the
  situation is trivial: $F$ is equal to the origin $\{0\}$ with a certain multiplicity $\omega$.
	But this multiplicity can be computed as $\omega = d_F(\RR^n)$. Hence, assuming that $d_F$ is the zero map,
	$\omega$ is zero as well. 
	
  To prove the induction step, we first use lemma
  \ref{SimplicialCompletion}, which
  shows that we can assume that $F$ can be represented as the $d$-skeleton
  of a complete simplicial rational fan $\Theta$ with certain (possibly zero)
  weights on the $d$-dimensional cones. We have to show that the
	assumption $d_F = 0$ implies that all these weights are zero.
  Let $\sigma$ be a $d$-cone of $\Theta$. As $\Theta$ is simplicial,
	we can find primitive vectors $v_1, \ldots, v_d$ that generate
  $\sigma$ 
	and a piecewise linear function $\varphi$ on $\Theta$ such that
	for each ray of $\Theta$ with primitive generator $v$ we have
	\begin{equation*} 
		\varphi(v) = \begin{cases}
		                 a \neq 0 & \text{for } v = v_1, \\
										 0 & \text{otherwise}.
		               \end{cases}
	\end{equation*}
	Let us now consider $\varphi \cdot F$. The compatibility
	of the divisor construction with the intersection product,
	i.e.\ $(\varphi \cdot F) \cdot Z = F \cdot (\varphi \cdot Z)$
	for all $Z \in Z_{n-d+1}(\RR^n)$, shows that 
	$d_{\varphi \cdot F} = 0$. 
	We apply the induction hypothesis and conclude that $\varphi \cdot F = 0$. 
	In particular, the weight $\omega_{\varphi \cdot F}(\tau)$ of  
  $\tau := \langle v_2, \ldots, v_d \rangle_{\RR_{\geq 0}}$ has to be zero.
	So let us compute this weight by hand:	
	Note that the primitive generator $v_{\sigma / \tau}$ of the projection of $\sigma$
	in $\RR^n / V_\tau$ is equal to the projection of $\frac{1}{|\Lambda_\sigma /
  (\Lambda_\tau + \ZZ v_1)|} v_1$ 
	(even though this vector itself might not be integer).
	Recall that $\varphi$ is identically zero on all facets
  containing $\tau$ except for $\sigma$ (in particular, $\varphi$
  is identically zero on $\tau$). Hence formula \eqref{weightdivisor} for the weight of $\tau$ 
	gives
	\begin{equation*} 
    \omega_{\varphi \cdot C}(\tau) =
      \omega_C(\sigma)
        \frac{1}{|\Lambda_\sigma / (\Lambda_\tau + \ZZ v_1)|}
        \varphi(v_1).
	\end{equation*}
	Since $|\Lambda_\sigma / (\Lambda_\tau + \ZZ v_1)|$ and $\varphi(v_1)$ are non-zero 
	numbers, $\omega_C(\sigma)$ must be zero, which finishes the proof.
\end{proof}

\section{The recession cycle}

Our goal is to compute the bounded Chow group $A^\text{b}_*(\RR^n)$ of $\RR^n$.
In proposition \ref{EquivalentCyclesAreEqual} we showed
that the group of fan cycles embeds into the bounded Chow group.
We will now show that the bounded Chow group is in fact isomorphic to the
group of fan cycles. To do so, we have to show that any tropical cycle
is bounded rationally equivalent to a fan cycle. Let us first describe this
(necessarily unique) fan cycle explicitly. 

\begin{definition}
  \label{def:recessioncycle}
  Let $\sigma$ be a polyhedron in $\RR^n$. We define the
  \dfn{recession cone} of $\sigma$ to be
	\begin{equation*} 
		\begin{split}
    \Rec(\sigma) :=& \; \{v \in \RR^n : x + \RR_{\geq 0} v \subseteq \sigma \; \forall \; x \in \sigma\} \\
                 =& \; \{v \in \RR^n : \exists \; x \in \sigma \text{ s.t. }
                                    x + \RR_{\geq} v \subseteq \sigma\}.
		\end{split}
	\end{equation*}
  The two sets coincide as $\sigma$ is closed and convex. 
  Let $X$ be a tropical $d$-dimensional cycle.
	It admits a polyhedral structure such $\XXX$ such that
	\begin{equation*} 
		\Rec(\XXX) := \{\Rec(\sigma) : \sigma \in \XXX\}
	\end{equation*}
	forms a fan, i.e.\ no cones overlap (cf.\ \cite[\nopp 1.4.10]{Ra2}). 
	We equip the $d$-cones of $\Rec(\XXX)$ with weights
	by 
	\begin{equation*} 
    \omega_{\Rec(\XXX)}(\sigma) :=
      \sum_{\substack{\sigma' \in X \\ \sigma = \Rec(\sigma')}}
      \omega_\XXX(\sigma').
	\end{equation*}
	This makes $\Rec(\XXX)$ a balanced fan (cf.\ \cite[61]{Ra2}) and we denote the
	corresponding cycle by $\Rec(X)$. We call $\Rec(\XXX)$ the \dfn{recession fan}
	of $\XXX$ and $\Rec(X)$ the \dfn{recession (fan) cycle} of $X$. Note that
	\begin{equation} \label{eq:recessionofsum}
	  \Rec(X + Y) = \Rec(X) + \Rec(Y).	
	\end{equation}
\end{definition}

\begin{example}
  \label{ex:recessionoftranslatedfan}
  Let $F$ be a fan cycle in $\RR^n$ and let $\vec{v} \in \RR^n$ be a vector. Then obviously
	\begin{equation*} 
		\Rec(F + \vec{v}) = F.
	\end{equation*}
	Indeed, when $\FFF = \{\sigma_i\}_i$ is a fan representing $F$, then $\FFF + \vec{v} = \{\sigma_i + \vec{v}\}_i$ 
	is a polyhedral structure for $F + \vec{v}$ and $\Rec(\FFF + \vec{v}) = \FFF$.
\end{example}

Our main result is the following:

\begin{theorem} \label{DegreeEquivalence}
  Let $X$ be a cycle in $\RR^n$. Then 
	\begin{equation*} 
		X \bsim \Rec(X).
	\end{equation*}
\end{theorem}

To prove this, we will use another theorem of its own interest. 

\begin{theorem}
  \label{thm:sumoffans}
  Let $X \in \RR^n$ be a tropical cycle. Then $X$ can be decomposed into 
	a sum of translated fan cycles, i.e.\ there are fan cycles $F_1, \ldots, F_l$
	and points $\vec{p}_1, \ldots, \vec{p}_l \in \RR^n$ such that
	\begin{equation*}
    X = \sum_{i=1}^l F_i + \vec{p}_i.
  \end{equation*}
\end{theorem}

The proof of this theorem (as it does not rely on the concept of rational equivalence)
will be postponed until section \ref{sec:decomposecycles}.
Instead, we continue with the proof of theorem \ref{DegreeEquivalence}, which of course 
is straightforward now.

\begin{proof}[Proof of \ref{DegreeEquivalence}]
  We write $X$ as a sum of translated fans $X = \sum_{i=1}^l F_i + \vec{p}_i$ by theorem \ref{thm:sumoffans}.
	By equation \ref{eq:recessionofsum} and example \ref{ex:recessionoftranslatedfan} we have
	\begin{equation*} 
		\Rec(X) = \sum_{i=1}^l F_i.
	\end{equation*}
	On the other hand, each translated fan $F_i + \vec{p}_i$ is bounded rationally equivalent to $F_i$ by
	proposition \ref{TranslationEquivalent}. As rational equivalence is additive, $X \bsim \Rec(X)$ 
	follows.
\end{proof}

Let us also mention another consequence of theorem \ref{thm:sumoffans}.
 
 \begin{corollary} \label{generatedbyhypersurfaces}
  Let $Z_*(\RR^n)$ denote the ring of tropical cycles in $\RR^n$, with $+$ the usual sum of cycles and $\cdot$ the stable intersection. Then $Z_*(\RR^n)$ is generated by 
	the set of hypersurfaces $V(f) \in Z_{n-1}(\RR^n)$ of tropical Laurent polynomials $f \in \TT[x_1^\pm, \ldots, x_n^\pm]$.
	In particular, $A_*(\RR^n) = 0$. 
  \begin{proof}
	  By \cite[\nopp 2.5.10]{MR}, every codimension one cycle can be written as a difference of two hypersurfaces $V(f) - V(g)$. Hence it suffices
		to show that $Z_*(\RR^n)$ is generated in codimension one. In the case of fan cycles, we can deduce this from the corresponding
		statement for (smooth) toric varieties and the equivalence of stable intersection and the toric intersection product
		(cf.\ \cite{FS, Ka, Ra}). Alternatively, a proof in purely combinatorial terms can be obtained via the polytope algebra
		(cf.\ \cite{FS,JYStable}).
		Finally, via theorem \ref{thm:sumoffans} we can reduce our case to the case of fan cycles and hence are done.
  \end{proof}
 \end{corollary}
 
 \begin{remark}
  In \cite{FS, JYStable}, the authors establish a link between the algebra of tropical fan cycles and 
	McMullen's polytope algebra \cite{MPolytope}.
	In the context of general cycles, one	can consider a generalized polytope algebra 
	generated by all polyhedra with a fixed given recession cone $\sigma$ 
	(in the ordinary case, $\sigma = \{0\}$). 
	Technically, this algebra might be constructed as a quotient of
	the ordinary polytope algebra by the additional relation
	\[
	  [P] = [Q] \textnormal{ if } P + \sigma = Q + \sigma.
	\]
	The case of interest for us are polytopes in $\RR^{n+1}$ and $\sigma =  \RR_{\geq 0} e_{n+1}$
	and hence the generators correspond, in some sense, to convex subdivisions of polytopes in $\RR^n$.
	Geometrically, this corresponds to taking tropical fan cycles in $\RR^{n+1}$ 
	and intersecting them with the hyperplane $\{x_{n+1} = - 1\}$. 
	Conjecturally, this generalized polytope algebra is isomorphic 
	to the algebra of general tropical cycles $Z_*(\RR^n)$ (not just fan cycles).
	However, the exact definitions and a subsequent proof of isomorphy still require careful analysis, 
	we do not pursue this here.

 \end{remark}

We finish this section by listing some consequences of theorem \ref{DegreeEquivalence}.
First, we conclude that the notions of bounded rational equivalence, 
numerical equivalence and ``having the same
recession cycle'' coincide.

\begin{theorem} \label{ComparisonOfRationalNumericalDegree}
  Let $X,Y$ be two tropical cycles in $\RR^n$. Then the
  following are equivalent:
  \begin{enumerate}
    \item[i)]
      $X \bsim Y$
    \item[ii)]
      $d_X = d_Y$
    \item[iii)]
      $\Rec(X) = \Rec(Y)$
  \end{enumerate}
	In particular, the equation
	\begin{equation*} 
		A^\text{b}_*(\RR^n) \cong Z_*^\text{fan}(\RR^n)
	\end{equation*}
	holds, where $A^\text{b}_*(\RR^n)$ is the bounded Chow group of $\RR^n$ and 
	$Z_*^\text{fan}(\RR^n)$ is the group of fan cycles.
\end{theorem}

\begin{proof}
  i) $\Rightarrow$ ii) follows from proposition
  \ref{RationalEquivalenceProperties} (d) and (e).
  iii) $\Rightarrow$ i) is an immediate consequence
  of theorem \ref{DegreeEquivalence}. ii) $\Rightarrow$ iii) follows
  from theorem \ref{DegreeEquivalence}, i) $\Rightarrow$ ii) and proposition
  \ref{EquivalentCyclesAreEqual}.
\end{proof}

The second corollary is the following general Bézout-type statement, where
$\Rec(X)$ plays the role of the degree of $X$.

\begin{theorem}[General Bézout's theorem]
  Let $X,Y$ be two tropical cycles in $\RR^n$. Then
  $$
    \Rec(X \cdot Y) = \Rec(X) \cdot \Rec(Y).
  $$
\begin{proof}
  We apply theorem \ref{DegreeEquivalence}
  and get
	\begin{equation*} 
    \Rec(X \cdot Y) \bsim X \cdot Y \bsim
    \Rec(X) \cdot \Rec(Y)
	\end{equation*}
  (the second equivalence also uses lemma
  \ref{RationalEquivalenceProperties} (d)).
  By proposition \ref{EquivalentCyclesAreEqual}
  two rationally equivalent fan cycles are equal.
\end{proof}
\end{theorem}

\section{Lineality spaces and splitting dimension}

In this section, we collect some additional definitions and notations which we need to prove
theorem \ref{thm:sumoffans}.

Let $X$ be a tropical cycle. A function $f : X \to \RR$ is called \dfn{lower semiconstant} if
for any polyhedral structure on $X$
\begin{enumerate}
	\item 
	  $f$ is constant on each relatively open cell $\RelInt(\sigma)$ (and hence we can 
		set $f(\sigma) := f(p)$, where $p$ is some point in the relative interior of $\sigma$),
	\item 
	  for any face $\tau \subseteq \sigma$ we have $f(\tau) \leq f(\sigma)$ (i.e.\
		$f$ is lower semicontinuous in the Euclidean topology). 
\end{enumerate}
Given such a function $f$ and $k \in \RR$, the sublevel set 
\begin{equation*}
  X_k = \{x \in X : f(x) \leq k\}
\end{equation*}
is again a polyhedral set.

Let $X$ be a tropical cycle and $p \in X$ a point. Locally around $p$,
$X$ looks like a fan and this fan cycle is denoted by $\Star_X(p)$.
As a set, $\Star_X(p)$ is the set of vectors $v \in \RR^n$ such that
$p + \epsilon v \in X$ for arbitrarily small $\epsilon > 0$. 
Given a polyhedral structure on $X$, we get an induced polyhedral structure
on $\Star_X(p)$ such that the facets of $\Star_X(p)$ are in one-to-one correspondence
to the facets of $X$ which contain $p$. Using the weights from $X$ for $\Star_X(p)$ accordingly,
the balancing condition is obviously still satisfied. 
Hence $\Star_X(p)$ is a fan cycle. It is easy to check the following formulas.
\begin{equation}
  \label{eq:starofsum}
  \Star_{X+Y}(p) = \Star_X(p) + \Star_Y(p)
\end{equation}
\begin{equation}
  \label{eq:starofstar}
  \Star_X(p+\epsilon q) = \Star_{\Star_X(p)}(q)
\end{equation}

Let $F \subseteq \RR^n$ be a fan cycle. The \dfn{lineality space} of $F$ is defined to be
\begin{equation*}
  \LinSp(F) := \{v \in \RR^n : F = F + v\},
\end{equation*}
where $F + v$ denotes the tropical cycle translated by $v$. 
Obviously, $\LinSp(F)$ is a linear subspace of $\RR^n$. Its dimension is denoted by
$\lindim(F)$ and is called the \dfn{lineality dimension} of $F$.
Examples are given in Figure \ref{linealityandsplittingdim} below.
In the special case $F = 0$ we set $\lindim(F) = \infty$.
When taking stars, we have
\begin{equation}
  \label{eq:linspacestars}
	\LinSp(F) \subseteq \LinSp(\Star_F(p))
\end{equation}
for all $p$.
When taking the sum of two fan cycles $F$ and $G$, we have
\begin{equation}
  \label{eq:linspacesum}
  \LinSp(F+G) \supseteq \LinSp(F) \cap \LinSp(G).
\end{equation}
We denote by
\begin{align*}
  \Fans^k &:= \{ F \text{ fan cycle} : \lindim(F) = k \}, \\
	\Fans^{\geq k} &:= \{ F \text{ fan cycle} : \lindim(F) \geq k \},
\end{align*}
the sets of fan cycles in $\RR^n$ with lineality space of dimension (greater than) $k$. 

\begin{definition}
  Let $F \subseteq \RR^n$ be a fan cycle. We define the \dfn{splitting dimension} of $F$ by
  \begin{equation*}
    \spldim(F) := \max\{k : F = \sum_i F_i \text{ for } F_i \in \Fans^{\geq k}\}.
  \end{equation*}
  Thus $\spldim(F)$ is the largest integer $k$ such that $F$ can be split into 
  a sum of fan cycles with lineality dimension at least $k$. 
  When $F = 0$, we have $\spldim(F) = \infty$.
  Let $X \subseteq \RR^n$ be a tropical cycle and let $p \in X$ be a point. We define
  the \dfn{lineality dimension} resp.\ \dfn{splitting dimension} of $p$ in $X$ by
  \begin{align*}
    l(p) &:= l_X(p) := \lindim(\Star_X(p)), \\ 
  	s(p) &:= s_X(p) := \spldim(\Star_X(p)). 
  \end{align*}
  In accordance with the previous conventions we set $l(p) = s(p) = \infty$ if $p \notin X$.
\end{definition}

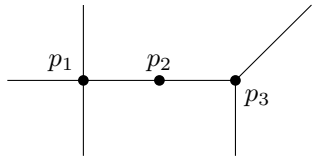
\begin{figure}[ht] 
  \centering
  \begin{tikzpicture}
    \draw (0,0) -- (3,0);
    \draw (1,1) -- (1,-1);
    \draw (3,0) -- (4,1);
    \draw (3,0) -- (3,-1);
    \fill[black] (1,0) circle (2pt) node[above left] {$p_1$};
    \fill[black] (2,0) circle (2pt) node[above] {$p_2$};
    \fill[black] (3,0) circle (2pt) node[below right] {$p_3$};
  \end{tikzpicture}
  \begin{minipage}{0.85\textwidth} 
    \caption{We calculate the lineality and splitting dimension of this one-dimensional tropical cycle at various points: 
	           $l(p_1) = l(p_3) = 0$, while $l(p_2) = 1$. $s(p_1) = s(p_2) = 1$, while $s(p_3) = 0$.}
  \label{linealityandsplittingdim}
	\end{minipage}
\end{figure}

Obviously, the chain of inequalities 
\begin{equation*}
  l(p) \leq s(p) \leq \dim(\Star_X(p))
\end{equation*}
holds and $l$ resp.\ $s$ are lower semiconstant functions on $X$ 
(by equation \eqref{eq:linspacestars}). 
It follows that the sets 
\begin{align*}
  X^{(k)} &:= \{ p \in X : l(p) \leq k\}, \\
	X^{[k]} &:= \{ p \in X : s(p) \leq k\},
\end{align*}
are polyhedral sets. We call $X^{(k)}$ the \dfn{$k$-skeleton} of $X$. 
Given a polyhedral structure for $X$, for each cell $\sigma$ we have
\begin{equation*}
  \dim(\sigma) \leq l(\sigma) \leq s(\sigma) 
\end{equation*}
and it follows
\begin{equation*}
  X^{[k]} \subseteq X^{(k)} \subseteq \bigcup_{\dim(\sigma) = k} \sigma. 
\end{equation*}
Moreover, these subsets are compatible with taking stars.

\begin{lemma}
  \label{lem:starskeletons}
	For any cycle $X$ and $p \in X$ we have
  \begin{align*}
    \Star_X(p)^{(k)} &= \Star_{X^{(k)}}(p), \\
	  \Star_X(p)^{[k]} &= \Star_{X^{[k]}}(p).
  \end{align*}
\end{lemma}

\begin{proof}
  Using equation \eqref{eq:starofstar}, we get the following chain of 
	equivalences.
  \begin{equation*}
	  \begin{split}
      q \in \Star_X(p)^{[k]} \; &\Longleftrightarrow \; 
		    \spldim(\Star_{\Star_X(p)}(q)) = \spldim(\Star_X(p+\epsilon q)) \leq k \\
			& \Longleftrightarrow \; p + \epsilon q \in X^{[k]} \\
			& \Longleftrightarrow \; q \in \Star_{X^{[k]}}(p)
		\end{split}
  \end{equation*}
	The case $\lindim$ is analogous.
\end{proof}

Here is another straightforward fact about lineality dimensions.

\begin{lemma}
  \label{lem:linsp-lindim}
  Let $F \subseteq \RR^n$ be a fan cycle and $p \in F$ a point. Then $l_F(p) \geq \lindim(F)$ and
	the equivalences
	\begin{equation*}
    l_F(p) = \lindim(F) \; \Longleftrightarrow \; p \in \LinSp(F) \; \Longleftrightarrow \; \Star_F(p) = F
  \end{equation*}
	hold.
\end{lemma}

\begin{proof}
  The inequality $l_F(p) \geq \lindim(F)$ is clear ($l_F$ is lower semiconstant). For
	the equivalences, we reduce to the case $\lindim(F) = 0$ by taking the quotient
	$F/\LinSp(F)$. Then the statement boils down to show
	\begin{equation*}
    l_F(p) = 0 \; \Longrightarrow \; p = 0 \; \Longrightarrow \; \Star_F(p) = F \; \Longrightarrow \; l_F(p) = 0.
  \end{equation*}	
	The first conclusion follows from the fact that each non-zero point in $F$ is contained
	in a positive-dimensional cell and therefore has positive lineality dimension. The remaining arrows are clear.
\end{proof}

\section{Decompose cycles into sums of fan cycles} \label{sec:decomposecycles}

In this section we prove theorem \ref{thm:sumoffans}, i.e.\ we show that every tropical cycle
can be decomposed into a sum of (translated) fan cycles.
The strategy of the proof is as follows. We recursively remove points in $X$ of minimal splitting
dimension by subtracting the corresponding star fans. The main step is to show that this
subtraction process does not create new points of minimal splitting dimension somewhere else. 
Based on this, we show that the process terminates (i.e.\ we obtain the zero-cycle) after a finite number of steps.

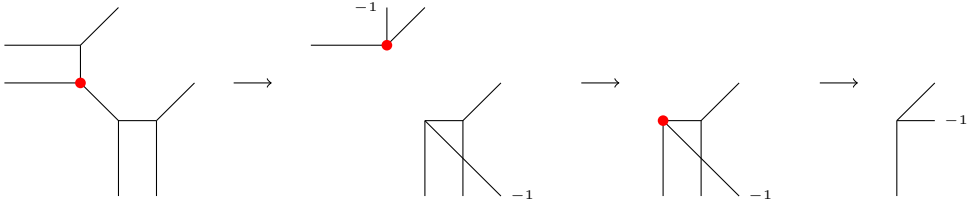
\begin{figure}[ht]
 \centering
 \begin{tikzpicture}[x=0.5cm, y=0.5cm, column sep = 0.5cm]
  \matrix{
      \draw (0,0) -- (0,-1) -- (1,-2) -- (2,-2);
      \draw (0,0) -- (-2,0);
      \draw (0,0) -- (1,1);
      \draw (0,-1) -- (-2,-1);
      \draw (1,-2) -- (1,-4);
      \draw (2,-2) -- (2,-4);
      \draw (2,-2) -- (3,-1);
      \fill[red] (0,-1) circle (2pt); &
      \draw[->] (0,-1) -- (1,-1); &
      \draw (0,0) -- (1,1);
      \draw (0,0) -- (-2,0);
      \draw (0,0) -- (0,1) node[left]{\tiny $-1$};
      \draw (1,-2) -- (2,-2) -- (3,-1);
      \draw (1,-2) -- (1,-4);
      \draw (2,-2) -- (2,-4);
      \draw (1,-2) -- (3,-4) node[right]{\tiny $-1$};
      \fill[red] (0,0) circle (2pt); &
      \draw[->] (0,-1) -- (1,-1); &
      \draw (1,-2) -- (2,-2) -- (3,-1);
      \draw (1,-2) -- (1,-4);
      \draw (2,-2) -- (2,-4);
      \draw (1,-2) -- (3,-4) node[right]{\tiny $-1$};
      \fill[red] (1,-2) circle (2pt); &
      \draw[->] (0,-1) -- (1,-1); &
      \draw (2,-2) -- (2,-4);
      \draw (2,-2) -- (3,-1);
      \draw (2,-2) -- (3,-2) node[right]{\tiny $-1$};
      \\
  };
 \end{tikzpicture}
  \begin{minipage}{0.91\textwidth}
    \caption{We iteratively subtract the local fans at the marked points from this plane tropical one-cycle until we arrive at a fan cycle. 
		         Note that in the first reduction step we obtain a new vertex on the bottom right leg of the one-cycle. 
						 This reduces the lineality dimension of this point to 0, but its splitting dimension is still 1.
						 (Of course, in the case of curves our approach is unnecessarily complicated and we could instead give an
						explicit formula for the decomposition into fans.)}
  \end{minipage}
\end{figure}

Let $X \subseteq \RR^n$ be a tropical cycle of dimension $\dim(X) = m$. We set
\begin{equation*}
  \sss := \sss(X) := \min\{s_X(p) : p \in \RR^n\}.
\end{equation*}
Our goal is to find a finite process which increases $\sss(X)$ by subtracting star fans.
It stops when $\sss(X) = \infty$ which means $X = 0$. The main step is contained in the 
following proposition.

\begin{proposition}
  \label{pro:unionofvectorspaces}
	The set of points of minimal splitting dimension $X^{[\sss]}$ is a finite union of affine subspaces $W_i \subseteq \RR^n$
	of dimension $\sss$,
  \begin{equation*}
    X^{[\sss]} = \bigcup_{i=1}^l W_i.
  \end{equation*}
\end{proposition}

To prove this, we use the following local condition.

\begin{lemma}
  \label{lem:localunionofvectorspaces}
  Let $Y \subseteq \RR^n$ be a polyhedral set such that for any $\vec{p} \in Y$ we have
  \begin{equation*}
    \Star_Y(\vec{p}) = \bigcup_{i=1}^{l_{\vec{p}}} V_{\vec{p},i},
  \end{equation*}
	where $V_{\vec{p},i} \subseteq \RR^n$ are $k$-dimensional linear subspaces. 
	Then $Y$ is a finite union of affine subspaces of dimension $k$.
\end{lemma}

\begin{proof}
  We choose a polyhedral structure for $Y$ with $k$-cells $\sigma_1, \ldots, \sigma_r$. 
	We pick a point $\vec{p}_j$ in the relative interior of $\sigma_j$ for all cells.
	Then $\Star_Y(\vec{p}_j) = V_j$ is a $k$-dimensional linear subspace.
	We want to show
  \begin{equation*}
    Y = \bigcup_{j=1}^r  V_j + \vec{p}_j,
  \end{equation*}
	i.e.\ $Y$ is equal to the union of translated spaces appearing in $\Star_Y(\vec{p}_j)$.
	The direction “$\subseteq$” is obvious as each cell $\sigma_j$ is contained in 
	$V_j + \vec{p}_j$. For the other inclusion, pick a vector space $V := V_j$ at a point
	$\vec{p}:= \vec{p}_j$. In order to show $V + \vec{p} \subseteq Y$, consider 
  \begin{equation*}
    Z := V \cap (Y-\vec{p}) \subseteq V.
  \end{equation*}
	By our assumptions, $Z$ is a full-dimensional polyhedral subset of $V$ with 
	$0 \in Z^\circ$. Assuming $Z \neq V$, we find a point $\vec{q}$ in the boundary of $Z^\circ$.
	But $\Star_Y(\vec{p}+\vec{q})$ is a union of $k$-dimensional vector spaces, and one of them must be $V$ as
	$\Star_Y(\vec{p}+\vec{q}) \cap V$ is $k$-dimensional. It follows $\vec{q} \in Z^\circ$, a contradiction. 
\end{proof} 

\begin{proof}[Proof of proposition \ref{pro:unionofvectorspaces}]
  Let $\sss := \sss(X)$ as above and pick $\vec{p} \in X^{[\sss]}$. By lemma \ref{lem:localunionofvectorspaces}
	it suffices to show that $\Star_{X^{[\sss]}}(\vec{p})$ is a union of linear subspaces of dimension $\sss$.
	
	We start by splitting $\Star_X(\vec{p})$ into a sums of fan cycles with lineality dimension at least $\sss$, or more
	precisely,
  \begin{equation}
	  \label{eq:startingsum}
    \Star_X(\vec{p}) = F_1 + \ldots F_l + \sum_i G_i,
  \end{equation}
	where $F_i \in \Fans^\sss$ and $G_i \in \Fans^{\geq \sss + 1}$.
	We set $V_i := \LinSp(F_i)$, hence $V_1, \ldots, V_l$ is a collection of 
	$\sss$-dimensional linear subspaces of $\RR^n$. We want to show that
	$\Star_{X^{[\sss]}}(\vec{p})$ is equal to a union of some of those $V_i$.
	First we show
  \begin{equation}
	  \label{eq:inclusioninunion}
    \Star_{X^{[\sss]}}(\vec{p}) \subseteq L := V_1 \cup \ldots \cup V_l.
  \end{equation}
	For all $q \in \RR^n$ we have
  \begin{equation}
	  \label{eq:starsum}
    \Star_{\Star_X(\vec{p})}(\vec{q}) = \sum_i \Star_{F_i}(\vec{q}) + \sum_i \Star_{G_i}(\vec{q})
  \end{equation}
	by equation \eqref{eq:starofsum}. If $\vec{q} \notin L$, then $\vec{q}$ is not contained
	in the lineality space of $F_i$ for all $i$ and thus 
	$\lindim(\Star_{F_i}(\vec{q})) > \sss$ for all $i$ (cf.\ lemma \ref{lem:linsp-lindim}).
	Hence on the right side of equation \eqref{eq:starsum}, all fans have 
	lineality dimension at least $\sss + 1$ and thus $s_{\Star_X(\vec{p})}(\vec{q}) \geq \sss +1$.
	We conclude $q \notin \Star_X(\vec{p})^{[\sss]}$, which by lemma \ref{lem:starskeletons}
	is the same as $\vec{q} \notin \Star_{X^{[\sss]}}(\vec{p})$. Equation \eqref{eq:inclusioninunion} follows.
	
	We now show the following: 
  \begin{equation}
	  \label{eq:notsubset}
    V_l \nsubseteq \Star_{X^{[\sss]}}(\vec{p}) \; \Longrightarrow \; \Star_{X^{[\sss]}}(\vec{p}) \subseteq V_1 \cup \ldots \cup V_{l-1}.
  \end{equation}
	This finishes the proof, as it allows us to recursively remove from equation \eqref{eq:inclusioninunion} all 
	vector spaces $V_i$ which are not contained in $\Star_{X^{[\sss]}}(\vec{p})$ until we reach equality. 
	To prove equation \eqref{eq:notsubset}, pick a point $\vec{q} \in V_l \setminus \Star_{X^{[\sss]}}(\vec{p})$. 
	Let us reorder the spaces $V_i$ (and $F_i$) such that 
	\begin{align*}
    \vec{q} &\notin V_i \text{ for all } i = 1, \ldots, r \\
	  \vec{q} &\in V_i \text{ for all } i = r+1, \ldots, l.
  \end{align*}
  We will prove the somewhat stronger statement $\Star_{X^{[\sss]}}(\vec{p}) \subseteq V_1 \cup \ldots \cup V_r$.
	Again by lemma \ref{lem:starskeletons}, we conclude from 
	$q \notin \Star_{X^{[\sss]}}(\vec{p}) = \Star_X(\vec{p})^{[\sss]}$ that $s_{\Star_X(\vec{p})}(\vec{q}) > \sss$. Thus we can write
  \begin{equation}
	  \label{eq:starstarsum}
    \Star_{\Star_X(\vec{p})}(\vec{q}) = \sum_i H_i
  \end{equation}
	for suitable fan cycles $H_i \in \Fans^{\geq \sss + 1}$.
	Combining equations \eqref{eq:starsum} and \eqref{eq:starstarsum} we get the expression
  \begin{equation}
	  \label{eq:replacementforFi}
    F_{r+1} + \ldots + F_l = \sum_i H_i - \left(\Star_{F_1}(\vec{q}) + \ldots + \Star_{F_r}(\vec{q}) + \sum_i \Star_{G_i}(\vec{q})\right).
  \end{equation}
	Here we used the fact that $\Star_{F_i}(\vec{q}) = F_i$ for all $i = r+1,\ldots,l$ by lemma \ref{lem:linsp-lindim}.
	Using this lemma again, we see that on the right hand side of this equation all fans have lineality dimension
	at least $\sss +1$. Finally, replacing the summands $F_{r+1} + \ldots + F_l$ in equation \eqref{eq:startingsum}
	by this expression, we get a new splitting of $\Star_X(p)$ of the form
  \begin{equation*}
    \Star_X(\vec{p}) = F_1 + \ldots F_r + \sum_i G'_i,
  \end{equation*}
	with fan cycles $G'_i \in \Fans^{\geq \sss + 1}$. Now the same reasoning as above 
	(which proved equation \eqref{eq:inclusioninunion}) shows that 
	$\Star_{X^{[\sss]}}(\vec{p}) \subseteq V_1 \cup \ldots \cup V_r$.
	This finishes the proof.
\end{proof}

Based on proposition \ref{pro:unionofvectorspaces} we now consider the
process of subtracting the star of a point of minimal splitting dimension.

\begin{proposition}
  \label{prop:subtractastar}
	Let $X \subseteq \RR^n$ be a tropical cycle with minimal splitting dimension $\sss := \sss(X)$.
	Write $X^{[\sss]} = W_1 \cup \ldots \cup W_l$, $W_i$ $\sss$-dimensional affine subspaces.
	Then there exists a point $\vec{p} \in W_l$ such that $W_l = \LinSp(\Star_X(\vec{p})) + \vec{p}$. 
	Moreover, for the tropical cycle $\widetilde{X} := X - (\Star_X(\vec{p}) + \vec{p})$, we have
  \begin{equation*}
    \widetilde{X}^{[\sss]} \subseteq W_1 \cup \ldots \cup W_{l-1}.
  \end{equation*}
\end{proposition}

\begin{proof}
  First, we show the existence of such a point $\vec{p}$. 
  Fix a polyhedral structure of $X$ and let $\sigma$ be a $\sss$-dimensional
	cell which is contained in $W_l$. Pick a point $p$ in the relative interior
	of $\sigma$. This implies $l_X(\vec{p}) \geq \dim(\sigma) = \sss$. But we also
	have $l_X(\vec{p}) \leq s_X(\vec{p}) = \sss$ and thus $l_X(\vec{p})=\sss$.
	Hence the lineality space of $\Star_X(\vec{p})$ is $\sss$-dimensional and is contained
	in $\Star_X(p)^{[\sss]} = \Star_{X^{[\sss]}}(\vec{p})$ (by lemma \ref{lem:starskeletons}).
	But $X^{[\sss]} = W_1 \cup \ldots \cup W_l$ and $p \in W_l$, hence $\Star_{X^{[\sss]}}(\vec{p}) = W_l - \vec{p}$, and we are done.

  Now let us check the second statement:
  Pick $\vec{q} \notin X^{[\sss]}$. From $\vec{q} \notin W_l$ it follows $\vec{q}-\vec{p} \notin \LinSp(\Star_X(\vec{p}))$ by assumption
	and	thus $\lindim(\Star_{\Star_X(\vec{p}) + \vec{p}}(\vec{q})) > \sss$ by lemma \ref{lem:linsp-lindim}.
	Write 
  \begin{equation*}
    \Star_X(\vec{q}) = \sum_i F_i
  \end{equation*}
	with fan cycles $F_i \in \Fans^{\geq \sss + 1}$. Using equation \eqref{eq:starofsum} we get
	\begin{equation*}
    \Star_{\widetilde{X}}(\vec{q}) = \sum_i F_i - \Star_{\Star_X(\vec{p}) + \vec{p}}(\vec{q}),
  \end{equation*}
	which implies $s_{\widetilde{X}}(\vec{q}) \geq \sss +1$ and $\vec{q} \notin \widetilde{X}^{[\sss]}$.
	This proves $\widetilde{X}^{[\sss]} \subseteq X^{[\sss]}$.
	By proposition \ref{pro:unionofvectorspaces} $\widetilde{X}^{[\sss]}$ must be equal
	to the union of some subcollection of the affine spaces $W_i$. Hence it suffices to show 
	$W_l \nsubseteq \widetilde{X}^{[\sss]}$. This follows from the fact that 
	by construction we have $\vec{p} \notin \widetilde{X}$,
	since $X$ and $\Star_X(\vec{p}) + \vec{p}$ coincide in a neighbourhood of $\vec{p}$. 
\end{proof}

We can now prove theorem \ref{thm:sumoffans}.

\begin{proof}[Proof of theorem \ref{thm:sumoffans}]
  We repeatedly subtract star fans as in proposition \ref{prop:subtractastar} in order to
	remove all points of splitting dimension $\sss$, i.e.\ $\widetilde{X}^{[\sss]} = \emptyset$.
	This implies $\sss(\widetilde{X}) > \sss$, and we can repeat the process until we reach
	$\sss(\widetilde{X}) = \infty$ and hence $\widetilde{X} = 0$ (alternatively, one may stop when
	$\sss(\widetilde{X}) = \dim(\widetilde{X})$ --- in this case $\widetilde{X} = \widetilde{X}^\sss$
	is a union/sum of affine subspaces). As during this procedure we only subtract translated
	fan cycles (namely of the form $\Star_X(\vec{p}) + \vec{p}$), the statement follows. 
\end{proof}

\printbibliography

\section*{Contact}

\begin{itemize}
	\item 
  Lars Allermann, Fachbereich Mathematik, TU Kaiserslautern, 
  Postfach 3049, 67653 Kaiserslautern, Germany; \href{mailto:allerman@mathematik.uni-kl.de}{allermanATmathematik.uni-kl.de}.

  \item
  Simon Hampe, Fachrichtung Mathematik, Universität der Saarlandes, 
  Postfach 151150, 66041 Saarbrücken, Germany; \href{mailto:hampe@math.uni-sb.de}{hampeATmath.uni-sb.de}.

  \item
  Johannes Rau, Fachrichtung Mathematik, Universität der Saarlandes, 
  Postfach 151150, 66041 Saarbrücken, Germany; \href{mailto:johannes.rau@math.uni-sb.de}{johannes.rauATmath.uni-sb.de}.
\end{itemize}

\end {document}